# Mathematics Graduate Training in the Age of AI

*Courses, Qualifying Exams, Oral Exams, and Dissertations*

David Glickenstein, University of Arizona

## Executive Summary

Generative AI changes the conditions under which graduate mathematics is learned, assessed, written, and defended. It can help students find examples, clarify definitions, generate code, improve exposition, and explore possible proof strategies. It can also produce confident but false mathematics, conceal gaps in understanding, and make it harder for conventional assignments to distinguish genuine mathematical command from tool-assisted answer production.

The central claim of this paper is that mathematics graduate programs should respond to the moment by clarifying what graduate mathematics education is trying to teach and assess. In most ways, the goals of mathematics education have not changed. Rather, with changing tools it has become more essential than ever to make clear the purpose of mathematical training.

Graduate mathematical learning has at least two dimensions. The first is local content mastery: definitions, theorems, examples, techniques, computations, and canonical proofs in a particular course or research area. The second is disciplinary mathematical formation: the portable practices by which mathematicians define, formalize, conjecture, prove, validate, communicate, and choose tools. While these elements of mathematical training are not new, AI makes the second dimension more visible and the need to clarify and cultivate it more urgent.

We use the term **mathematical judgment** to refer to the capacity to evaluate mathematics (e.g., claims, definitions, examples, proofs, analogies, computations, uses of tools, research directions) as mathematically sound, useful, well-posed, and appropriately justified. The phrase is not introduced here as a wholly new term: related uses appear in the mathematics education literature on mathematical competency, mathematical thinking, reasoning, problem handling, and more, e.g., Niss and Højgaard (2019). In this paper we use the phrase as an organizing construct for graduate assessment in an AI-rich environment.

The literature on mathematical competencies, proof comprehension, expert/novice proof reading, reading mathematical exposition, problem solving, metacognition, and advanced mathematical thinking gives language for assessing mathematical maturity. The literature on AI literacy, prompt engineering, overreliance on AI, human-AI collaboration, and AI assessment design gives usable language for assessing expert versus non-expert AI use. Taken together, these literatures suggest that students should be assessed not only on what they can produce, but on more: what they can read, monitor, verify, repair, explain, transfer, and defend.

The practical recommendation is to adjust gradually: keep tests, qualifying exams, oral exams, courses, proposals, and dissertation defenses, but fold judgment tasks into each format. Written exams can include flawed-proof diagnosis, missing-hypothesis repair, definition comparison, counterexample construction, and proof-structure explanation. Oral exams can probe why a definition is natural, where a hypothesis is used, what breaks without it, and what a student would try next. Dissertations can require explicit defense of definitions, examples, proof architecture, and any AI-assisted work.

Again, these are not new ideas; the recommendation is simply to emphasize these techniques and make their use explicit.

AI policy should be differentiated by academic setting. Graduate courses may allow AI for explanation, examples, coding, writing feedback, and exploration when students disclose use and verify output. Timed qualifying exams should normally be AI-free unless explicitly redesigned as AI-permitted assessments. Dissertation research and writing should require committee-approved use, disclosure, and independent verification of any AI-assisted material. Across all settings, the student remains responsible for the mathematics submitted under their name. Recent mathematics-specific guidance emphasizes the same core principles: transparency, human responsibility, verification, and preservation of disciplinary judgment (Leiden 2026).

## Acknowledgements

This document originally arose out conversations on training incoming Math PhD students, primarily with Adi Adiredja. It was greatly influenced by conversations and correspondences on training students and on AI with many people, including Bryden Cais, Misha Chertkov, Ben Chow, Yi Hu, Mac Hyman, Joceline Lega, Kevin Lin, Feng Luo, Selena Niu, Marek Rychlik, Patrick Shipman, Hang Xue, and Doug Ulmer.

## Disclosure of AI Use in Generating This Document

AI was used in generating this document. Specifically, an AI language model was used to synthesize the prior discussion, organize the argument, draft prose, suggest assessment structures, and assemble a preliminary bibliography and policy landscape. The author has verified all citations and policy descriptions, reviewed and revised all of the prose, and takes full responsibility for the final claims and recommendations.

## 1. Introduction: AI and the Changing Conditions of Graduate Mathematics Education

Mathematics graduate programs have long relied on a relatively stable set of educational pieces, including graduate courses, problem sets, written qualifying examinations, oral examinations, reading courses, research seminars, dissertation proposals, dissertations, and defenses. These elements of a graduate program are designed to train and evaluate students to become mathematicians: people who can formulate precise questions, use definitions well, build and assess proofs, choose examples/nonexamples/counterexamples, recognize structure, and participate in the research community.

Generative AI stresses this ecology. Large language models can produce polished mathematical prose, plausible proof outlines, code, examples, bibliographic suggestions, and explanations at a speed that changes ordinary student work. Students can use AI to avoid work they do not wish to do. In doing so, AI gives the appearance of mathematical fluency, independent of correctness or clarity of exposition. For instance, a proof can look coherent while containing a false lemma, a missing hypothesis, a fabricated citation, or a subtle quantifier error. A dissertation paragraph can sound polished while weakening the logical structure of an argument. An important point is that many of these issues are not unique to AI, but could be done by a person; however, the ability of AI to produce fluent but incorrect mathematics *at scale* makes the issue more pressing and apparent.

This paper tackles the question of how graduate mathematics programs can use AI responsibly as a productivity and learning tool. In graduate mathematics, AI should be governed by the same standards that govern all serious mathematical work. As emphasized in Leiden (2026), mathematics is a human enterprise depending on expert judgment of mathematicians, and should be attributable to and verified by human authors. From the perspective of education, a student may use tools, but the student must take responsibility for the mathematics.

The guiding question is: What should mathematics graduate programs assess in an AI-rich environment? The answer proposed here is that programs should explicitly assess mathematical judgment. Judgment includes the capacity to determine whether a claim is well-posed, whether a proof proves what it says it proves, whether a definition captures the intended structure, whether an example is relevant, whether an analogy is sound, whether a computational result is trustworthy, and whether an AI-generated response should be accepted, revised, or rejected.

# 2. What Graduate Mathematics Programs Currently Assess

Current graduate assessment practices are not uniform, but they tend to emphasize several recognizable forms of evidence. Graduate courses use homework, exams, projects, and sometimes oral presentations to evaluate command of material. Written qualifying examinations test fluency with core definitions, theorems, computations, and proof methods. Oral examinations test recall, explanation, adaptability, and live reasoning. Dissertation proposals and defenses test research direction, independence, and ownership of mathematical work.

These formats have important strengths. They require students to engage with core material, demonstrate technical proficiency, and communicate complicated concepts. They also have limitations. A written qualifying exam can reward rehearsed proof production without revealing whether a student can diagnose a flawed proof. Homework can reward polished final solutions without revealing false starts, strategy choice, or testing who did the actual work.

The problem with many of these formats is that their most important goals are implicit. Faculty often say that a student is or is not mathematically mature, but the criteria for that judgment are rarely stated with the same precision as course syllabi or exam topics. Braun discusses how the term mathematical maturity is common but imprecise, involving technical skill, habits of investigation, persistence, conceptual understanding, and broader intellectual functioning (Braun, 2019).

The AI era makes the weakness of such vague criteria even more critical. If students can obtain plausible solutions, proof outlines, and exposition from AI systems, then assessment must be clearer about what is being evaluated: the answer, the method, the verification, the explanation, the judgment, or the student's independent command.

# 3. Mathematical Maturity, Mathematical Practices, and Mathematical Judgment

We make a distinction here between local content mastery and disciplinary mathematical formation. Local content mastery concerns the material at hand: the definitions, theorems, examples, counterexamples, standard constructions, computations, and proof techniques of a course or research area. Disciplinary mathematical formation concerns the portable practices by which mathematicians work across fields: defining, formalizing, conjecturing, proving, validating, reading, connecting, communicating, and choosing tools.

In graduate mathematics these practices are essential to learning and communicating. The ability to turn text into quantifiers, to see which hypothesis is doing real work, to read a proof for its architecture, or to recognize a structure across fields is part of mathematical knowledge itself.

The KOM framework (Competencies and the Learning of Mathematics, Niss and Højgaard, 2019) attempts to separate mathematical competence, or *competency*, from particular content areas. Niss and Højgaard describe mathematical competency in terms of activities such as mathematical thinking, problem handling, modeling, reasoning, representation, use of symbols and formalism, communication, and use of aids and tools (Niss and Højgaard, 2019). For graduate programs, core math classes, such as algebra or real analysis, share many program-level learning goals despite covering different content.

The proof-comprehension literature gives another attempt at clarifying mathematical reasoning. Mejia-Ramos and collaborators propose a multidimensional model of proof comprehension that includes the meaning of terms and statements, the logical status of claims, the chaining of implications, the high-level ideas, the proof method, the modular structure, and the relation of the proof to examples (Mejia-Ramos et al., 2012). Selden and Selden's work on proof validation and Inglis and Alcock's comparison of expert and novice proof readers similarly show that knowing what a proof says is not the same as judging whether it proves the theorem (Selden and Selden, 2003; Inglis and Alcock, 2012).

The reading of mathematical exposition is also a disciplinary practice. Shepherd and van de Sande (2014) compare first-year undergraduates, graduate students, and faculty mathematicians reading unfamiliar mathematical exposition and propose a Mathematics Reading Framework with three dimensions: mathematical fluency, comprehension monitoring, and engagement. In their study, faculty mathematicians differed from graduate students and undergraduates by more often reading the meaning of symbols, monitoring comprehension, persevering through confusing passages, using figures and other resources, and constructing examples while reading. This difference is not solely due to more content knowledge. Difference can occur within a single individual as well. Students may be advanced in one area while remaining novice or intermediate readers in another.

Schoenfeld's account of mathematical problem solving adds another dimension: students need resources, strategies, metacognitive monitoring, and productive beliefs about mathematics (Schoenfeld, 1992/2016). This matters especially in graduate study, where open-ended work requires choosing what to try, noticing when an approach is failing, and revising a plan. Tall's edited volume explores the jump to advanced mathematical thinking, summarized as:

*The move from elementary to advanced mathematical thinking involves a significant transition: that from describing to defining, from convincing to proving in a logical manner based on those definitions. This transition requires a cognitive reconstruction (Tall, 1991, p. 20).*

Metacognition gives a complementary vocabulary for this development. Tanner (2012) emphasizes that students can be taught to plan, monitor, and evaluate their learning, rather than leaving "learning how to learn" to chance. The MIT Teaching + Learning Lab similarly describes metacognition as using knowledge of the task, strategies, and oneself to plan learning, monitor progress, and evaluate outcomes (MIT Teaching + Learning Lab, n.d.). Medina, Castleberry, and Persky (2017) connect metacognition to critical thinking, self-regulated learning, questioning, modeling, reflection, and judgments of understanding. Persky and Robinson (2017) frame expertise development as staged and dependent on progressive problem solving, organized knowledge, and self-regulation. Boaler et. al. (2024) suggest ways that metacognition can be a powerful tool both teaching/learning and assessing mathematics. These accounts support treating mathematical knowledge and learning as partly technical and partly metacognitive: students need to know mathematics, but they also need to monitor how they are reading, proving, computing, and deciding.

For purposes of AI-era learning and assessment, we consider the unifying term mathematical judgment. Mathematical judgment is the capacity to evaluate mathematical objects and arguments as mathematics: to decide whether a definition is well-shaped, whether a claim is true, whether a proof is valid, whether an example is relevant, whether an analogy is legitimate, whether a computation supports the conclusion, and whether a tool output deserves trust.

Some of the learning goals and assessment possibilities for mathematical judgment are listed in Table 1. Table 1 is a synthesis rather than a reproduction from a single source. The categories draw on the KOM framework for mathematical competencies (Niss and Højgaard, 2019), models of proof comprehension and validation (Mejia-Ramos et al., 2012; Selden and Selden, 2003; Inglis and Alcock, 2012), Schoenfeld's work on problem solving and metacognition (Schoenfeld, 1992/2016), research on reading mathematical exposition from novice to expert (Shepherd and van de Sande, 2014), metacognition frameworks centered on planning, monitoring, and evaluating (Tanner, 2012; Medina et al., 2017; MIT Teaching + Learning Lab, n.d.), and AI-literacy work on critical, accountable tool use (Long and Magerko, 2020; Passi and Vorvoreanu, 2022).

### Table 1. Learning goals and possible evidence

| Learning goal | What the student should be able to do | Possible assessment evidence |
|---|---|---|
| Local content mastery | State definitions, use standard examples, prove central theorems, and solve representative problems. | Problem sets, exams, theorem-proof questions, computations, and oral explanation. |
| Formalization | Turn prose into hypotheses, conclusions, quantifiers, diagrams, and symbolic statements. | Formalize-this-statement tasks; identify missing hypotheses; compare proposed formalizations. |
| Proof validation | Decide whether a proof is valid and locate the first serious gap or false step. | Flawed-proof diagnosis, proof annotation, live proof repair, AI proof audit. |
| Proof comprehension | Explain the main idea, dependency structure, key lemma, and proof method. | One-paragraph proof summaries, dependency diagrams, oral proof navigation. |
| Mathematical reading | Read definitions, notation, examples, and exposition for meaning; monitor comprehension; use examples, diagrams, and resources to resolve confusion. | Guided reading annotations, read-the-meaning tasks, notation-to-prose explanations, example construction while reading, reading-to-learn interviews. |
| Definition sense | Judge whether a definition is natural, too weak, too strong, circular, or equivalent to another definition. | Definition comparison, example/non-example/counterexample tests, definition-design assignments. |
| Transfer and synthesis | Recognize patterns across fields while respecting limits of analogy. | Cross-topic problems, synthesis essays, seminar presentations, comparative oral questions. |
| Metacognitive regulation | Plan an approach, monitor progress and confusion, evaluate outcomes, and revise strategies when evidence shows that an approach is not working. | Exam wrappers, learning journals, reflective memos, confidence judgments with explanation, plan-monitor-evaluate checklists. |
| Tool-mediated judgment | Use AI, computation, search, or proof assistants while retaining responsibility for verification. | AI-use appendix, reproducible code, verification checklist, defense of tool-generated claims. |
| Research formation | Formulate questions, judge significance and feasibility, locate work in literature, and defend choices. | Literature map/annotated bibliography, dissertation proposal, defense, referee-style report. |

# 4. Expert Versus Non-Expert Use of AI

The distinction between expert and non-expert AI use is the distinction between calibrated, accountable use and trusting but ignorant dependence. Expert use is selective, critical, transparent, and verified. Non-expert use accepts plausible output too quickly, loses track of sources, fails to test edge cases, and mistakes fluency for correctness.

The novice-to-expert literature helps describe this distinction. Expertise is domain-specific, staged, and tied to organized knowledge, progressive problem solving, and metacognitive regulation (Persky and Robinson, 2017). Expert use of AI in mathematics depends on the user's ability to notice what is mathematically relevant, recognize when a tool response does not fit the definitions or hypotheses, and decide what must be checked independently.

Long and Magerko define AI literacy through competencies that enable people to critically evaluate, communicate, and collaborate with AI systems (Long and Magerko, 2020). Chiu and collaborators distinguish AI literacy from AI competency, emphasizing not only knowledge about AI but also the ability to apply that knowledge with confidence and self-reflection (Chiu et al., 2024). Recent work on generative AI literacy assessment, including the GLAT instrument, reinforces the need for performance-based measures rather than self-report alone (Jin et al., 2025).

AI produces answers through prompting, and thus prompt engineering is an important piece of working with AI. However, prompt engineering is only part of AI literacy. Federiakin and collaborators frame prompt engineering as the skill of articulating a problem, its context, and constraints to an AI assistant (Federiakin et al., 2024). Knoth and collaborators connect AI literacy to prompt-engineering strategies (Knoth et al., 2024). However, even though literacy produces better prompts, literacy also involves evaluation of the resulting output. This is crucial for the use of AI in mathematics.

Hyman (2026) makes the important distinction between the tool-use skill and the learning with AI skill. To actually learn with AI requires knowing when and how to challenge AI answers, think about the concepts oneself, and follow up.

The literature on overreliance on AI looks at how the user handles the output of AI. Microsoft's Aether review defines overreliance as accepting incorrect AI outputs and emphasizes the need for effective explanations and mitigation of automation bias (Passi and Vorvoreanu, 2022). Klingbeil and collaborators provide experimental evidence that trust and reliance can lead users to follow AI advice even when it is costly or inappropriate, underscoring the need for trust calibration (Klingbeil et al., 2024).

Generative AI can also give uneven help. Brynjolfsson, Li, and Raymond find that generative AI can increase productivity in customer support, with especially large gains for novice and lower-skilled workers (Brynjolfsson et al., 2023/2025). Dell'Acqua and collaborators describe a 'jagged technological frontier' in which AI helps on some tasks but worsens performance on others, even within the same professional workflow (Dell'Acqua et al., 2023). The implication for mathematics is that a graduate program should not assume either that AI is uniformly helpful or uniformly harmful. It should teach students to recognize where the frontier lies for a given mathematical task.

Work comparing experts and novices suggests caution about treating expertise as a simple binary. Sun and collaborators (2022) found that novices and domain experts contributed differently to chatbot-improvement tasks: novices could perform comparably on some lower-level classification tasks, while experts contributed more to tasks requiring domain judgment, such as creating new intent categories and authoring contextually appropriate responses. For mathematics graduate programs, the analogous point is that some AI-supported tasks may be made safer through verification routines, while others require domain expertise and disciplinary judgment.

For graduate mathematics, expert AI use can be characterized by behaviors: precise problem framing; awareness of AI limitations; independent verification; use of examples, nonexamples, and counterexamples; transparent disclosure; and final ownership of the submitted mathematics. These behaviors can be assessed directly.

## 5. Existing AI Policies in Mathematics and Graduate Education

We review publicly available AI policies to illustrate several workable governance models. It is important to note that such policies change frequently. The Leiden Declaration (2026) gives an overview of values for mathematics research that can serve as a general guidepost for Mathematics, drawing heavily on the UNESCO Recommendation on Open Science (2021).

The University of Notre Dame Department of Mathematics has a policy that applies to graduate comprehensive exams, theses, doctoral scholarly projects, and dissertations. Its model is close to a prohibition-with-exceptions approach: AI-generated content is prohibited in these high-stakes settings unless the use is approved under specified conditions, such as committee approval or accommodation, and approved use must be disclosed (Notre Dame Department of Mathematics, n.d.).

Oxford's Mathematical Institute provides a discipline-specific model. Its MSc dissertation guidance and departmental policy distinguish permissible uses such as literature search, code assistance, formatting, grammar, and ideation from impermissible uses such as substantive AI-generated writing, undeclared AI-generated code, production of plots that obscure data or algorithms, and direct use of AI-generated interpretations of mathematics or data (Oxford Mathematical Institute, 2024/2025). The policy recognizes that mathematical writing is assessed not merely for grammatical correctness but for clarity and logical structure.

KIT's Department of Mathematics has also issued AI guidelines for mathematics students and teachers as a supplement to institution-wide guidance (KIT Department of Mathematics, 2026). Northwestern's mathematics graduate syllabus statements provide a course-level model: generative AI may be used as supplementary material to aid understanding, but homework must be written in the student's own words and checked by the student; asking an LLM for a full answer is discouraged (Northwestern Department of Mathematics, n.d.).

Graduate-school policies at Georgia, Auburn, Illinois State, and Arizona State provide committee-centered models. Georgia treats generative AI use in theses and dissertations as unauthorized unless specifically authorized by the advisory committee within an approved scope and requires disclosure if approved. Auburn requires AI use in graduate research and writing to be disclosed and approved by the advisory committee. Illinois State requires disclosure of all AI use in thesis, dissertation, capstone, and comprehensive exam processes to the student's committee. Arizona State places responsibility on the graduate student's committee to determine permissible amounts and types of AI use in culminating experiences (University of Georgia Graduate School, 2024; Auburn University Graduate School, n.d.; Illinois State Graduate School, n.d.; Arizona State University Graduate College, n.d.).

There are also some general higher-education frameworks available. The AI Assessment Scale offers a practical way for instructors to specify the level of AI use permitted in a given assessment, from no AI through AI-supported work and AI exploration (Perkins et al., 2024). TEQSA's work on assessment reform emphasizes that higher education should redesign assessment to manage AI risks while supporting responsible use (TEQSA, 2023/2025). University-level syllabus policy resources from institutions such as Ohio State, Florida, and Texas give adaptable language for instructors, though they usually need disciplinary refinement for graduate mathematics.

Table 2 gives a summary of some policy models based on the descriptions given in this section.

## Table 2. Possible policy models for mathematics graduate programs

| Policy model | Typical rule | Use in mathematics graduate programs |
|---|---|---|
| Prohibition with exceptions | AI use is prohibited in high-stakes work unless explicitly approved. | Appropriate for timed qualifying exams, comprehensive exams, and final dissertation text unless redesigned. |
| Committee-approved disclosure | Advisor or committee determines allowable use; student discloses use and scope. | Well-suited to dissertations, proposals, research writing, and individualized projects. |

| Discipline-specific allowed use | Policy distinguishes editing, search, code, plots, mathematical interpretation, proof generation, and authorship. | Acknowledges that different uses carry different risks for mathematical rigor and clarity. |
|---|---|---|
| Learning support but not final answer | AI may help explain, practice, or debug, but submitted work must be independently written and checked. | Useful for graduate courses and qualifying-exam preparation. |
| AI-integrated verification | AI use is permitted, but students must audit, verify, cite, and defend all accepted output. | Useful for assignments specifically designed to teach expert AI use and mathematical judgment. |

An important critique of many policy models is given by Corbin et. al. (2025), who make a distinction between discursive changes, ones that depend on the student to comply with rules, vs. structural changes, ones that change the nature of how the assessment is being made. Most policy changes give a careful description of what is allowed or not. However, without a change to how assessment is made it is unclear how valid the assessment is. The authors suggest that structural changes often require a shift from testing output to testing process. Also, it makes sense to have valid assessments at a less granular level, such as at the level of units or modules, rather than at a finer level such as individual assignments. This fits well into graduate program goals of *training mathematicians* rather than completing a homework set or individual qualifying exam. In the following sections, we should keep in mind this distinction between discursive and structural changes to assessment. Notice that most of policies in Table 2 are of the discursive type, indicating a need to consider more structural solutions.

# 6. AI in Graduate Courses

Graduate courses are a good place to teach responsible AI in a formative way. Students can experiment, make mistakes, receive feedback, and learn habits of verification before reaching high-stakes milestones. A course policy should specify what kinds of uses for AI are allowed and what evidence of student understanding is expected and tested.

In a graduate course, AI may reasonably be permitted for asking for alternative explanations, generating practice problems, suggesting examples, improving grammar, helping with code syntax, exploring computations, or producing first-pass summaries of background material. These uses should not replace student responsibility for the final mathematical argument. A student should not submit an AI-generated proof as if it were their own work; nor should a student rely on an AI-generated theorem, reference, or computation without independent verification. It should be more like they received the work from a colleague as a collaboration.

A basic curricular change is to add judgment tasks to ordinary coursework. After a proof problem, ask students to identify where each hypothesis was used. After a theorem, ask for a counterexample when one hypothesis is removed. After a definition, ask for two nearby definitions and examples distinguishing them. After an AI-assisted exploration, ask what was accepted, what was rejected, and how the accepted material was checked.

Courses can also use AI output as mathematical text to be read, criticized, and repaired. This turns AI from an invisible shortcut into an object of disciplinary instruction. Students learn that a fluent proof is not necessarily a valid proof, that a plausible definition may be too weak or too strong, and that mathematical responsibility cannot be delegated to a tool.

Courses can also teach mathematical reading explicitly. A reading assignment might ask students to translate notation into prose, identify the role of each definition, mark where comprehension breaks down, construct a small example, and explain what resource they would use next. This work makes visible a practice that experienced mathematicians often perform silently.

Metacognitive routines can be folded into existing assignments with little structural change. A course can ask students before an assignment to name the task, resources, and likely obstacles; during the assignment to record confusions or failed approaches; and after feedback to explain what they would do differently. Exam or other metacognitive wrappers, confidence judgments, and short reflections are examples of such routines (Tanner, 2012; Medina et al., 2017; Lovett, 2023, MIT Teaching & Learning Lab).

Hyman (2026) suggest assignments may be given of different types, all of which have an important place in education of students:

- No AI. AI is not to be used at all. These are meant to test the ability of the student without AI. Examples include closed-book exams/quizzes, in-class writing, oral explanations, live coding checks, language drills, or foundational calculations.
- AI as tutor. AI can be used to learn but not to produce content. Permitted AI use includes explaining a concept, quiz the student, give hints, create practice problems, or help prepare for class.
- AI as collaborator. AI can be used with brainstorming, drafting, debugging, revision, critique, outlining, coding, or comparison. The results must still be verified and disclosed.
- AI as object of study. AI output is looked at and evaluated directly. Students may compare AI responses, find errors, test generated code, identify bias, repair a flawed explanation, or evaluate source work.

### Possible course-level assignments

- Flawed proof diagnosis: give a short proof with a subtle quantifier error, a false implication, or an unproved compactness claim. Ask students to find the first serious error, explain it, and repair the statement or proof.
- Definition comparison: give two definitions that are close but not equivalent. Ask students to decide whether they are equivalent and to produce examples distinguishing them.
- Formalization from prose: give an informal mathematical statement and ask for a precise version with quantifiers, hypotheses, and conclusion.
- Proof idea/proof detail: require both a three-sentence structural explanation and a complete proof.
- AI proof audit: ask an AI system for a proof, then require students to annotate the response, identify unsupported steps, and replace them with valid arguments from course material.
- AI-use appendix: for projects, require a brief statement of what AI was used for, what was independently checked, and what was rejected.
- Mathematical reading audit: assign a short passage from a text or paper and ask students to translate notation, identify definitions and dependencies, construct examples, and list unresolved questions.
- Metacognitive wrapper: after an exam, proof assignment, or project, ask students what they expected to know, where their understanding failed, what evidence they used to evaluate their performance, and what they will change next time.

## 7. AI and Qualifying Exams

Qualifying examinations serve several purposes: they certify foundational knowledge, mark progress toward candidacy, and often operate as gatekeeping mechanisms. Recent work emphasizes that qualifying exams are widespread but under-studied and that evidence for standardized effective practices remains limited (McLaughlin et al., 2024). AI makes it especially important to state what qualifying exams are meant to certify.

Timed written qualifying exams should be AI-free unless explicitly redesigned with AI in mind. The purpose of such exams is to evaluate independent command under constraint. If AI is allowed without redesign, the exam may no longer measure the intended construct. But an AI-free exam should not mean a traditional exam only. Qualifying exams can remain closed-tool while still assessing mathematical judgment more directly.

One practical model is to reserve a portion of the exam, perhaps 20 to 30 percent, for judgment-oriented tasks. The rest can remain standard proof, computation, and problem solving. Judgment tasks might include flawed-proof diagnosis, theorem repair, example/counterexample construction, comparison of definitions, reading-to-learn from a short unfamiliar passage, and explanation of proof structure. These questions reveal aspects of graduate readiness that purely content-based questions may not show.

If a program chooses to experiment with explicitly AI-permitted qualifying components, the learning outcomes should be adjusted appropriately. An AI-permitted component ask for more than a final solution to a problem. It should ask students to use AI to generate candidate arguments or examples, then audit, verify, correct, and defend the result. The assessment target would be expert tool use and mathematical judgment and distinguish this from simple recall.

### Sample qualifying-exam judgment prompts

- Analysis: A student claims that pointwise convergence plus convergence of integrals implies uniform convergence. Evaluate the claim. If false, give a counterexample, then state a nearby true theorem.
- Topology: A proposed proof asserts that compactness is preserved under arbitrary intersections by passing finite subcovers across the family. Identify the flaw and give the corrected statement.
- Algebra: Two definitions of a normal subgroup are given. Determine whether they are equivalent, prove equivalence if they are, or give a distinguishing example if not.
- PDE: A formal energy estimate is given. Identify the regularity assumptions needed to justify each integration by parts step.
- Reading task: A short excerpt from an unfamiliar but accessible mathematical text is provided. Translate the notation into prose, identify the main definition and its hypotheses, construct a simple example, and state one question that would need to be resolved before using the result.

## 8. AI and Oral Exams

Oral examinations are particularly valuable in an AI-rich environment because they can assess live mathematical ownership. A student may prepare with books, notes, peers, code, or AI, but in the room they must explain, respond, repair, and connect. Oral exams therefore test proficiency in a way written exams do not.

Oral exams should include standard content questions, but they should also include explicit judgment prompts. After a definition, ask why that definition is natural and what breaks if a condition is removed. After a proof, ask for the main idea, the key lemma, the use of hypotheses, and a possible generalization. After a theorem, ask for examples at the boundary of the hypotheses. After a failed argument, ask the student to diagnose and repair it.

Oral exams can also assess mathematical reading in real time. The examiner can give a definition, theorem statement, or short proof excerpt and ask the student to read the notation as a mathematical sentence, identify dependencies, explain why a phrase is needed, or construct an example illustrating the statement. This is close to what students must do when reading papers, using AI-generated mathematical text, or entering a new research area.

Hyman (2026) suggests asking students to re-present work in another form, e.g., different notation, graphically, a plain-language explanation, or a different example. These are ways to test whether the student has mastered the underlying idea rather than just having recognized an already polished answer.

The oral format also allows examiners to probe AI-era issues without making the exam about AI. For example, an examiner can present a plausible but false proof and ask the student to evaluate it. This resembles the task a student will face when reading AI-generated mathematics, but it is also simply good mathematical training.

### Core oral-exam judgment questions

- Why is this the right definition?
- Where exactly is this hypothesis used?
- What example shows that the hypothesis cannot simply be removed?
- What is the main idea of the proof, separate from the technical details?
- Which lemma is the heart of the argument?
- Could this proof work in a more general setting? Why or why not?
- What analogy from another field is helpful here, and where does the analogy break?
- Suppose this approach fails. What would you try next?
- Can you read this line of notation as a mathematical sentence, and say what each symbol is doing?
- Illustrate this concept graphically.
- Here is a similar but different setting with different notation. Explain the same concept but with these notations.

## 9. AI in Dissertation Research and Writing

In dissertation work the distinction between tool use and authorship becomes especially important. Mathematicians routinely use tools: computation, databases, search engines, symbolic systems, proof assistants, typesetting systems, collaborators, seminars, conversations and, increasingly, AI. AI should be understood within this broad ecology of mathematical work, paying special attention to its capacity to generate plausible text, code, references, and arguments that may not be reliable or attributable in ordinary scholarly ways. This is consistent with the recommendations of the Leiden declaration (2026).

A mathematics dissertation policy should distinguish categories of AI use such as literature search, bibliographic formatting, grammar and style editing, code generation or debugging, plotting or visualization, example generation, conjecture or strategy generation, and proof or exposition generation. These categories should not be treated alike. Grammar checking and bibliography formatting are low-risk when disclosed and checked. AI-generated proof, mathematical interpretation, or substantive exposition is high-risk because it may alter the mathematical content itself.

One governance model for dissertation work is committee-approved disclosure. The student and advisor should discuss proposed AI use early, ideally at the proposal stage. The committee should decide which uses are permitted, which require documentation, and which are prohibited. Any final dissertation should include a disclosure statement describing the role of AI tools in research, coding, writing, editing, or visualization.

All of the previous safeguards are of the discursive type as described by Corbin et. al. (2025). However, there is already a structural safeguard, the dissertation defense. The dissertation defense should assess mathematical ownership. If AI was used to explore examples, write code, search literature, improve exposition, or suggest proof strategies, the student should be able to explain what was used, what was checked, what was discarded, and why the final argument is correct. No mathematical claim belongs in the dissertation unless the student can defend it independently. It may become more important for the dissertation defense to significantly assess correctness and attribution.

### Suggested dissertation tool and computational resource use disclosure

The following paragraph may be adapted for dissertation proposals, dissertations, or project reports:

**Tool and Computational Resource Use Disclosure.** In preparing this work, I used generative AI tools for the following purposes: [list purposes, such as literature search, coding assistance, grammar editing, organization, example exploration, or drafting feedback]. I did not use AI-generated text, code, mathematical arguments, figures, or references in the final document without independent verification. All mathematical claims, proofs, computations, citations, and interpretations in the final document are my responsibility. [If applicable: Specific AI-assisted materials are identified in Appendix X, together with a description of how they were checked.]

## 10. Recommendations for Policy and Assessment Design

The following recommendations are intended for mathematics departments, graduate committees, and faculty designing courses, qualifying exams, oral exams, and dissertation policies. They are deliberately modular to allow different programs to adopt them according to size, culture, field distribution, and institutional constraints.

1. **Make mathematical judgment an explicit program outcome.** Graduate programs should state that students are expected not only to know content but also to evaluate definitions, proofs, examples, conjectures, computations, analogies, and tool outputs.
2. **Distinguish local content mastery from disciplinary mathematical formation.** Course syllabi and exam guidelines should name both the local material and the portable practices being assessed.
3. **Teach mathematical reading as a graduate skill**. Programs should ask students to read definitions, notation, examples, and proofs for meaning, not only to reproduce finished arguments.
4. **Build metacognitive routines into graduate courses**. Assignments can ask students to plan, monitor, and evaluate their work, especially around difficult reading, proof repair, and AI-assisted exploration.
5. **Add judgment tasks to existing exams rather than replacing existing exams.** Written tests and qualifying exams can include flawed proofs, missing hypotheses, formalization, proof summaries, example/counterexample tasks, and definition comparisons.
6. **Use oral exams to assess ownership and adaptability.** Oral exams should include questions about why definitions are natural, where hypotheses are used, how proofs work structurally, and what examples show the limits of a theorem.
7. **Create differentiated AI policies by academic setting.** Homework, projects, timed exams, oral exams, proposals, dissertations, and defenses should not all have the same AI rule.
8. **Require disclosure and verification for AI-assisted work.** Students should say how AI was used, what was accepted, what was rejected, and how any accepted output was checked.
9. **Do not rely on AI detection as the main integrity mechanism.** Assessment should be designed around secure evidence of learning, oral defense, process documentation, and verification rather than detection of AI-generated text.
10. **Treat AI-generated mathematics as untrusted until verified.** A fluent proof, reference, computation, or explanation should be checked against definitions, accepted theorems, examples, source literature, or reproducible computation.
11. **Use AI as an object of instruction.** Courses can teach judgment by asking students to critique AI-generated proofs, definitions, examples, and explanations.
12. **Make dissertation committees responsible for field-sensitive AI guidance.** Because acceptable AI use varies by mathematical area and research method, dissertation committees should approve scope, documentation, and disclosure expectations early in the process.

## Conclusion

AI is becoming part of the environment in which students read, write, compute, and learn. This moment provides an opportunity to clarify and strengthen what and how we assess in graduate programs.

The central recommendation is to place mathematical judgment at the center of graduate education. Students should still learn definitions, theorems, techniques, and proofs. They should still solve problems without assistance. However, emphasis should be placed on learning to judge mathematics: to read critically, monitor their own understanding, verify carefully, define precisely, use examples strategically, understand proof structure, recognize the benefits and limits of analogy, and govern their tools responsibly.

In the end, a mathematics graduate program should aim to form mathematicians who can decide what is correct, explain why it is correct, recognize when it is not, and take responsibility for the mathematics they develop and present to others.

## Appendix A. A Rubric for Mathematical Judgment

Note that Mathematical Judgment as described here touches on all elements of Hyman's paradigm of core competence, AI-assisted practice, and trustworthiness.

| Dimension | Weak evidence | Strong evidence |
|---|---|---|
| Precision | Uses vague terms; omits hypotheses or quantifiers. | States claims with appropriate hypotheses, quantifiers, and conclusions. |
| Validation | Accepts plausible arguments too quickly. | Checks logical dependencies and identifies gaps or false steps. |
| Examples | Uses vague or irrelevant examples. | Uses examples, nonexamples, and counterexamples strategically to clarify ideas and to test claims. |
| Proof structure | Focuses only on local steps. | Identifies main idea, key lemma, and dependency structure. |
| Definition sense | Memorizes definitions without explaining their role. | Explains what the definition captures, for instance through examples/ nonexamples and analogy. |
| Mathematical reading | Reads notation verbatim, skips difficult passages, or cannot identify where comprehension failed. | Translates notation into meaning, checks comprehension, constructs examples or diagrams, and seeks targeted clarification. |
| Transfer | Treats topics as isolated. | Recognizes structures across fields while noting limits of analogy. |
| Metacognitive regulation | Equates time spent with learning; lacks a plan for monitoring or revising strategy. | Plans a strategy, monitors progress and confusion, evaluates results, and revises the approach based on evidence. |
| Repair | Can say that something is wrong but cannot explain why or fix it. | Adds hypotheses, corrects statements, or replaces invalid steps. |
| Tool judgment | Trusts AI or computation uncritically. | Verifies tool outputs and explains what was checked independently or still needs verification. |

## References and Policy Sources